\documentclass[12pt,leqno]{article}
\usepackage{amssymb,amsthm,amsmath,bbm}
\def\bm#1{\mathbbm{#1}}
\def\fn#1{\mathop{{\rm #1}\vphantom{\dim}}}

\makeatletter
\renewcommand{\section}{\@startsection{section}{1}{0mm}{12mm}{5mm}{\raggedright\bf\Large}}
\makeatother

\newtheorem{theorem}{\indent Theorem}[section]

\newtheorem{example}[theorem]{\indent Example}
\newtheorem{prop}[theorem]{\indent Proposition}

\title{Reduction theorems for Noether's problem}
\author{{ Ming-chang Kang}\\[2mm]
Department of Mathematics and\\Taida Institute of Mathematical Sciences\\ National Taiwan University\\
Taipei, Taiwan\\ E-mail: kang@math.ntu.edu.tw
\\[3mm] and
\\[3mm] { Bernat Plans}\\[2mm]
Departament de Matem\`{a}tica Aplicada I\\ Universitat
Polit\`{e}cnica de Catalunya\\ Barcelona, Spain\\ E-mail:
bernat.plans@upc.edu}

\date{}

\begin{document}

\maketitle
\footnote{\hspace*{-7.5mm} 2000 Mathematics Subject Classification. Primary 12F12, 12F20, 13A50, 11R32, 14E08.}

\begin{abstract}
\noindent Abstract. Let $K$ be any field, $G$ be a finite group.
Let $G$ act on the rational function field $K(x(g):g\in G)$ by
$K$-automorphisms and $h\cdot x(g)=x(hg)$. Denote by
$K(G)=K(x(g):g\in G)^G$ the fixed field. Noether's problem asks
whether $K(G)$ is rational (= purely transcendental) over $K$. We
will give several reduction theorems for solving Noether's
problem. For example, let $\widetilde{G}=G\times H$ be a direct
product of finite groups. Theorem. Assume that either (1) $H$ is
an abelian group of exponent $e$ and $K$ contains a primitive
$e$-th root of unity, or (2) $K$ is a field with $\fn{char} K=p >
0$ and $H$ is a $p$-group. Then $K(\widetilde{G})$ is rational
over $K(G)$. In particular, if $K(G)$ is rational (resp.\ retract
rational) over $K$, so is $K(\widetilde{G})$ over $K$.
\end{abstract}

\newpage
\section{Introduction}

Let $K$ be any field, $G$ be a finite group.
Let $G$ act on the rational function field $K(x(g):g\in G)$ by $K$-automorphisms and $h\cdot x(g)=x(hg)$.
Denote by $K(G)=K(x(g):g\in G)^G$ the fixed field.
Noether's problem asks whether $K(G)$ is rational (= purely transcendental) over $K$.
For a survey of Noether's problem, see Swan's paper \cite{Sw}.

The purpose of this article is to prove several reduction theorems
when we try to solve Noether's problem for some group. First we
will prove the following theorem without assuming Fischer's
Theorem (see Theorem \ref{t1.2}).

\begin{theorem} \label{t1.1}
Let $\widetilde{G}=H\times G$ be a direct product of finite
groups, and let $K$ be a field. Assume that {\rm (i)} $H$ is an
abelian group with exponent $e$, i.e. $e=\max\{\fn{ord}(h):h\in
H\}$; {\rm (ii)} the field $K$ contains a primitive $e$-th root of
unity. Then there is a $K$-embedding of $K(G)$ into
$K(\widetilde{G})$ so that $K(\widetilde{G})$ is rational over
$K(G)$.
\end{theorem}

By a $K$-embedding of $K(G)$ into $K(\widetilde{G})$ we mean an
injective $K$-linear homomorphism of fields from $K(G)$ into
$K(\widetilde{G})$. Note that, for any field $K$, if $G$ and
$\widetilde{G}$ are finite groups so that $K(\widetilde{G})$ is
rational over $K(G)$, then $K(\widetilde{G})$ is rational (resp.\
stably rational, retract rational) over $K$ provided that so is
$K(G)$. (Recall that ``rational" $\Rightarrow$ ``stably rational"
$\Rightarrow$ ``retract rational". For the definition of retract
rationality, see \cite [Definition 3.2]{Sa2}.) Thus Theorem
\ref{t1.1} becomes a very convenient technique in solving
Noether's problem or proving the existence of generic
$G$-polynomials. An immediate consequence of Theorem \ref{t1.1} is
the classical Fischer's Theorem.

\begin{theorem} \label{t1.2} {\bf (Fischer's Theorem \cite[Theorem 6.1]{Sw})}
Let $G$ be a finite abelian group of exponent $e$, and let $K$ be
a field containing a primitive $e$-th root of unity. Then $K(G)$
is rational over $K$.
\end{theorem}

A result similar to Theorem \ref{t1.1} when $\fn{char} K=2$ is the following.

\begin{theorem} \label{t1.3} {\bf (\cite[Proposition 7]{Pl})}
Let $K$ be a field with $\fn{char} K=2$ and $\widetilde{G}$ be a
group extension defined by $1\to \bm{Z}/2\bm{Z} \to
\widetilde{G}\to G\to 1$ where $G$ is a finite group. Then
$K(\widetilde{G})$ is rational over $K(G)$.
\end{theorem}

Combining Theorem \ref{t1.1} and Theorem \ref{t1.3}, we obtain the
following result.

\begin{theorem} \label{t1.4}
Let $K$ be any field, and $\widetilde{G}=(\bm{Z}/2\bm{Z})\times G$
be a direct product of finite groups. Then $K(\widetilde{G})$ is
rational over $K(G)$.
\end{theorem}

Another application of Theorem \ref{t1.1} is the case of dihedral
groups, for which we will denote by $D_n$ the dihedral group of
order $2n$. The following theorem is implicit in \cite{Ka}.

\begin{theorem} \label{t1.5}
If $K$ is any field and $n$ is an odd integer, then $K(D_{2n})$ is
rational over $K(D_n)$. In particular, if $K(D_n)$ is rational
{\rm (}resp.\ retract rational{\rm )} over $K$, so is $K(D_{2n})$.
\end{theorem}

\begin{proof}
If $D_{2n}=\langle \sigma,\tau: \sigma^{2n}=\tau^2=1,\
\tau\sigma\tau^{-1}=\sigma^{-1}\rangle$, then $D_{2n}$ is a direct
product of the groups $\langle \sigma^2,\tau \rangle$ and $\langle
\sigma^n \rangle$. Apply Theorem \ref{t1.4}. Note that $\langle
\sigma^2,\tau \rangle$ is isomorphic to $D_n$.
\end{proof}

Here is a generalization of Theorem \ref{t1.3} to the case when
$\fn{char} K=p$.

\begin{theorem} \label{t1.6}
Let $K$ be a field with $\fn{char} K=p>0$ and $\widetilde{G}$ be a
group extension defined by $1\to \bm{Z}/p\bm{Z} \to
\widetilde{G}\to G\to 1$ where $G$ is a finite group. Then
$K(\widetilde{G})$ is rational over $K(G)$.
\end{theorem}

An application of the above theorem is the following.

\begin{theorem} \label{t1.7}
Let $K$ be a field with $\fn{char} K=p>0$ and
$\widetilde{G}=H\times G$ be a direct product of finite groups
where $H$ is a $p$-group. Then there is a $K$-embedding of $K(G)$
into $K(\widetilde{G})$ so that $K(\widetilde{G})$ is rational
over $K(G)$.
\end{theorem}

\begin{proof}
Induction on the order of $H$. Let $\sigma \in H$ be an element of
order $p$ and $\sigma$ is contained in the center of $H$. Define
$G'=(H/<\sigma>) \times G$. Then we get a short exact sequence
$1\to <\sigma> \to \widetilde{G}\to G'\to 1$. Apply Theorem
\ref{t1.6}. We find that $K(\widetilde{G})$ is rational over
$K(G')$.
\end{proof}

A corollary of the above theorem is Kuniyoshi's Theorem : If $K$
is a field with $\fn{char} K=p > 0$ and $G$ is a finite $p$-group,
then $K(G)$ is rational over $K$ \cite{Ku}.

We record another application of Theorem \ref{t1.6}.

\begin{theorem} \label{t1.8}
Let $K$ be a field with $\fn{char} K=p>0$ and $\widetilde{G}$ be a
group extension defined by $1\to H \to \widetilde{G}\to G\to 1$
where $H$ and $G$ are finite groups. If $H$ is a cyclic $p$-group
or an abelian $p$-group lying in the center of $\widetilde{G}$,
then $K(\widetilde{G})$ is rational over $K(G)$.
\end{theorem}

Finally we will give two variants (or generalizations) of Theorem
\ref{t1.1}.
\begin{theorem} \label{t1.9}
Let $K$ be any field, and $H$ and $G$ be finite groups. If $K(H)$
is rational {\rm(}resp. stably rational, retract rational {\rm)}
over $K$, so is $K(H\times G)$ over $K(G)$.

In particular, if both $K(H)$ and $K(G)$ are rational {\rm(}resp.
stably rational, retract rational{\rm)} over $K$, so is $K(H\times
G)$ over $K$.
\end{theorem}

\begin{theorem} \label{t1.10}
Let $K$ be any field, $H\wr G$ be the wreath product of finite
groups $H$ and $G$. If $K(H)$ is rational {\rm(}resp. stably
rational {\rm)} over $K$, so is $K(H\wr G)$ over $K(G)$.
\end{theorem}

Note that it is known that, for an infinite field $K$, if $K(H)$
and $K(G)$ are retract rational over $K$, so are $K(H\times G)$
and $K(H\wr G)$ over $K$ (\cite[Theorem 1.5 and Theorem
3.3]{Sa1}). An application of Theorem \ref{t1.9} and Theorem
\ref{t1.10} to Noether's problem for dihedral groups will be given
in Theorem \ref{t4.2}.

We will prove Theorem \ref{t1.1}, Theorem \ref{t1.6}, Theorem
\ref{t1.9} and Theorem \ref{t1.10} in Section 2, Section 3, and
Section 4 respectively.

\medskip
Standing notations. We will denote by $\zeta_n$ a primitive $n$-th
root of unity. When we say that a field $K$ contains a primitive
$n$-th root of unity, it is assumed tacitly that $\fn{char} K=0$
or $\fn{char}K=p>0$ with $p\nmid n$. If $G$ is a finite group, we
will write $V=\bigoplus_{g\in G} K\cdot x(g)$ as the regular
representation space of $G$ where $G$ acts on $V$ by $h\cdot
x(g)=x(hg)$ for any $g,h\in G$. Recall the definition
$K(G)=K(x(g):g\in G)^G$ defined at the beginning of this section.

\bigskip

\bigskip

\section{Proof of Theorem \ref{t1.1}}

Before proving Theorem \ref{t1.1} we recall two basic facts.

\begin{theorem}\label{t2.1} {\bf (Hajja and Kang \cite[Theorem 1]{HK})}
Let $G$ be a finite group acting on $L(x_1,\ldots,x_n)$,
the rational function field of $n$ variables over a field $L$.
Suppose that

{\rm (i)} for any $\sigma \in G$, $\sigma (L)\subset L$;

{\rm (ii)} the restriction of the action of $G$ to $L$ is faithful;

{\rm (iii)} for any $\sigma \in G$,
\[
\left(\begin{array}{c} \sigma(x_1) \\ \sigma(x_2) \\ \vdots \\ \sigma(x_n) \end{array}\right)
=A(\sigma)\cdot \left(\begin{array}{c} x_1 \\ x_2 \\ \vdots \\ x_n \end{array}\right) +B(\sigma)
\]
where $A(\sigma)\in GL_n(L)$ and $B(\sigma)$ is an $n\times 1$ matrix over $L$.

Then there exist $z_1,\ldots,z_n\in L(x_1,\ldots,x_n)$ so that $L(x_1,\ldots,x_n)=L(z_1,\ldots,z_n)$
with $\sigma(z_i)=z_i$ for any $\sigma\in G$, any $1\le i \le n$.
\end{theorem}

\begin{theorem} \label{t2.2} {\bf (Ahmad, Hajja and Kang \cite[Theorem 3.1]{AHK})}
Let $L$ be any field, $L(x)$ the rational function field of one
variable over $L$, and $G$ a group acting on $L(x)$. Suppose that,
for any $\sigma\in G$, $\sigma(L)\subset L$, and
$\sigma(x)=a_\sigma \cdot x+b_\sigma$ where $a_\sigma,b_\sigma\in
L$ and $a_\sigma\ne 0$. Then $L(x)^G=L^G$ or $L^G(f)$ for some
polynomial $f\in L[x]$. In fact, if the integer $m:=\min\{\deg
g(x):g(x)\in L[x]^G,\ g(x)\notin L\}$ does exist, then
$L(x)^G=L^G(f(x))$ for any $f(x)\in L[x]^G$ satisfying $\deg f=m$.
\end{theorem}

\bigskip

\bigskip

\begin{proof}[\indent Proof of Theorem \ref{t1.1} ] ~

Step 1. Suppose that Theorem 1.1 is valid when $H$ is a cyclic
group. Then it is also valid when $H$ is an abelian group, because
we may write $H$ as a direct product of cyclic groups and use
induction on the number of these cyclic groups.

>From now on, we will assume that $H$ is a cyclic group of order
$n$.

\medskip

Step 2. Write $H=\langle c\rangle$ and $\zeta=\zeta_n$. Write the
coset decomposition $\widetilde{G}=\bigcup_{g\in G} gH$.

Let $\widetilde{V}=\bigoplus_{g\in \widetilde{G}} K\cdot
x(\tilde{g})$ and $V=\bigoplus_{g\in G} K\cdot x(g)$ be the
regular representation spaces of $\tilde{G}$ and $G$ respectively.

\medskip

Step 3. For each $g\in G$, define
\[
z(g)=\sum_{0\le i\le n-1} \zeta^i x(c^i g)\in \widetilde{V}.
\]

Define
\[
W=\bigoplus_{g\in G} K\cdot z(g) \subset \widetilde{V}.
\]

Note that $h\cdot z(g)=z(hg)$, $c\cdot z(g)=\zeta^{-1} z(g)$ for any $g,h\in G$.
It follows that $\widetilde{G}$ acts faithfully on $K(z(g):g\in G)$.
Apply Theorem \ref{t2.1} to $K(z(g):g\in G)$ and $K(x(\tilde{g}):\tilde{g}\in \widetilde{G})$.
We find that $K(\widetilde{G})$ is rational over $K(z(g):g\in G)^{\widetilde{G}}$.

\medskip

Step 4. If $G=\{1\}$, the trivial group, then $K(z(g):g\in
G)^{\widetilde{G}}=K(z(1)^n)$ is rational over $K$. From now on,
we assume that $G$ is not the trivial group.

\medskip

Step 5. For each $h\in G\backslash \{1\}$, define
\[
t(h)=z(h)/z(1).
\]

It follows that $K(z(g):g\in G)=K(t(h):h\in G\backslash \{1\})$
($z(1))=L(z(1))$ where $L=K(t(h):h\in G\backslash\{1\})$. Note
that, for any $g\in G$, $g\ne 1$,
\begin{eqnarray}
&& g\cdot z(1)=z(g)=(z(g)/z(1))z(1),\ c\cdot z(1)=\zeta^{-1}z(1), \label{eq2.1} \\
&& g\cdot t(h)=t(gh)/t(g)\in L,\ c\cdot t(h)=t(h). \label{eq2.2}
\end{eqnarray}

Because of (\ref{eq2.1}) and (\ref{eq2.2}), we may apply Theorem
\ref{t2.2}. Hence $K(z(g):g\in
G)^{\widetilde{G}}=L^{\widetilde{G}}(t_0)$ for some $t_0$ with
$\tilde{g}\cdot t_0=t_0$ for any $\tilde{g}\in \widetilde{G}$.

Because of (\ref{eq2.2}), we find that $L^{\widetilde{G}}=L^G$.
Thus $K(z(g):g\in G)^{\widetilde{G}}=K(t(h):h\in G\backslash \{1\})^G (t_0)$.

\medskip

Step 6. Consider $K(G)=K(x(g):g\in G)^G$. For each $h\in
G\backslash\{1\}$, define
\[
s(h)=x(h)/x(1).
\]

It follows that $K(x(g): g\in G)=K(s(h):h\in G\backslash
\{1\})(x(1))=L'(x(1))$ where $L'=K(s(h):h\in G\backslash\{1\})$.
Note that, for any $g\in G$, $g\ne 1$,
\begin{eqnarray}
&& g\cdot x(1)=(x(g)/x(1))\cdot x(1), \label{eq2.3} \\
&& g\cdot s(h)=s(gh)/s(g)\in L. \label{eq2.4}
\end{eqnarray}

Imitate the trick in Step 4. We find that $K(G)={L'}^G(s_0)$ for
some $s_0$ with $g\cdot s_0=s_0$ for any $g\in G$. Moreover,
$K(G)=K(s(h):h\in G\backslash\{1\})^G(s_0)$. Compare (\ref{eq2.2})
and (\ref{eq2.4}). We find that $K(t(h):h\in G\backslash
\{1\})^G(t_0)$ is $K$-isomorphic to $K(s(h):h\in
G\backslash\{1\})^G(s_0)$.
\end{proof}

\bigskip

\bigskip

\begin{example} \label{ex2.3} \em
The assumption that $\zeta_e \in K$ in Theorem \ref{t1.1} is
crucial.

In fact, let $\widetilde{G}=\bm{Z}/8\bm{Z} \times \bm{Z}/4\bm{Z}$
and $G=\bm{Z}/4\bm{Z}$. Then $\bm{Q}(G)$ is rational, but
$\bm{Q}(\widetilde{G})$ is not even retract rational \cite[Theorem
5.11]{Sa1}.
\end{example}

\begin{example} \label{ex2.4} \em
We don't know whether Theorem \ref{t1.1} is valid for
$\widetilde{G}$ which is a semi-direct product, but not a direct
product. In fact, we don't know whether there exist distinct prime
numbers $p$ and $q$ such that $\widetilde{G}=\bm{Z}/p\bm{Z}
\rtimes \bm{Z}/q\bm{Z}$ is a non-abelian semi-direct product and
$\bm{C}(\widetilde{G})$ is not rational over $\bm{C}$.

However, consider the non-abelian group
$\widetilde{G}=\bm{Z}/17\bm{Z} \rtimes \bm{Z}/16\bm{Z}$ where
$\bm{Z}/16\bm{Z}$ acts faithfully on $\bm{Z}/17\bm{Z}$. By Serre's
Theorem \cite[Theorem 33.16, p.88]{GMS}, $\bm{Q}(\widetilde{G})$
is not retract rational over $\bm{Q}$ (and neither is
$\bm{Q}(\bm{Z}/16\bm{Z})$ by \cite{Sa1}), while it is known that
both $\bm{C}(\widetilde{G})$ and $\bm{C}(\bm{Z}/16\bm{Z})$ are
rational over $\bm{C}$ \cite[Theorem 3.5]{Sa1}.
\end{example}

\begin{example} \label{ex2.5} \em
We may even try to work out a result similar to Theorem \ref{t1.1}
for the case of a non-split group extension in view of Theorem
\ref{t1.6}. But this is impossible. Just consider the extension
$0\to \bm{Z}/2\bm{Z} \to \bm{Z}/8\bm{Z} \to \bm{Z}/4\bm{Z} \to 0$.
Note that $\bm{Q}(\bm{Z}/4\bm{Z})$ is rational over $\bm{Q}$ while
$\bm{Q}(\bm{Z}/8\bm{Z})$ is not retract rational over $\bm{Q}$
\cite[Theorem 5.11]{Sa1}.
\end{example}

\bigskip

\bigskip

\section{Proof of Theorem \ref{t1.6}}

In this section, $K$ is a field with $\fn{char} K=p > 0$ and $1\to
\bm{Z}/p\bm{Z} \to \widetilde{G}\to G\to 1$. Let $c$ be a
generator of the normal subgroup $\bm{Z}/p\bm{Z}$ and $\pi
:\widetilde{G}\to G\to 1$ be the given epimorphism.

The idea of the proof is somewhat similar to the proof of Theorem
\ref{t1.1}.

\medskip
Step 1. Let $u:G \to \widetilde{G}$ be a section of $\pi$.

As before let $\widetilde{V}=\bigoplus_{g\in \widetilde{G}} K\cdot
x(\tilde{g})$ and $V=\bigoplus_{g\in G} K\cdot x(g)$ be the regular
representation spaces of $\tilde{G}$ and $G$ respectively.

\medskip
Step 2. For each $g\in G$, define
\[
y(g)=\sum_{0\le i\le p-1} x(c^i u(g))\in \widetilde{V},
\]

\[
z(g)=\sum_{0\le i\le p-1} i x(c^i u(g))\in \widetilde{V},
\]

\[
z=\sum_{g \in G}z(g)\in \widetilde{V},
\]

\[
W=(\bigoplus_{g\in G} K\cdot y(g)).
\]

Note that $c\cdot y(g)=y(g)$. As $G$-spaces, $W$ and $V$ are
$G$-equivariant. Hence $K(W)^{\widetilde{G}} \simeq K(G)$.

\medskip
Step 3. We will examine the action of $\widetilde{G}$ on $z(g)$ and
$z$.

It is clear that $c\cdot z(g)=z(g)-y(g)$.

For any $h,g \in G$, suppose that $u(h)\cdot u(g)=c^m\cdot u(hg)$
and $u(h)\cdot c \cdot u(h)^{-1}=c^n$. Note that $m$ is an integer
depending on $g$ and $h$, and $n$ is invertible in $K$. When the
element $h$ is fixed, we may write $m=m(g)$ to emphasize the
dependence of $m$ on $g$.

We find that $u(h)\cdot z(g)=\sum_{0 \le i \le p-1} i x(u(h)c^i
u(g))$ $= \sum_{0 \le i \le p-1} i x(c^{in}u(h) u(g))$ $= \sum_{0
\le i \le p-1} i x(c^{in+m}u(hg))$ $=c^m \cdot (1/n)\sum_{0 \le i
\le p-1} i x(c^{i}u(hg))$ $=(1/n) z(hg)-(m/n)y(hg).$

It follows that $u(h)\cdot z=(1/n) z -\sum_{g \in G}(m(g)/n )y(hg)$
where $m(g)$ denotes the integer $m$ depending on $g$.

\medskip
Step 4. Define $\widetilde{W}=W \bigoplus K \cdot z$. Then
$\widetilde{W}$ is a faithful $\widetilde{G}$-subspace of
$\widetilde{V}$. By Theorem \ref{t2.1}, $K(\widetilde{G})$ is
rational over $K(\widetilde{W})^{\widetilde{G}}$.

Consider the pair $\widetilde{W}$ and $W$ and apply Theorem
\ref{t2.2}. We find that $K(\widetilde{W})^{\widetilde{G}}$ is
rational over $K(W)^{\widetilde{G}}$. Since
$K(W)^{\widetilde{G}}=K(W)^G \simeq K(G)$, we are done.

\bigskip

\bigskip

\section{Proof of Theorem \ref{t1.9} and Theorem \ref{t1.10}}

\begin{proof}[\indent Proof of Theorem \ref{t1.9}] ~

Without loss of generality, we may assume that neither $H$ nor $G$
is the trivial group.

Step 1. Write $\tilde{G}=H\times G$.

Let $U=\bigoplus_{h\in H} K\cdot x(h)$ and $V=\bigoplus_{g\in G}
K\cdot x(g)$ be the regular representation spaces of $H$ and $G$
respectively.

For any element $\tilde{g} \in \widetilde{G}$, any $u \otimes v
\in U \bigotimes_K V$, define $\tilde{g} \cdot (u \otimes v)=(h
\cdot u)\otimes (g \cdot v)$ if $\tilde{g}=hg$ where $h \in H$ and
$g \in G$. It is easy to see that $U \bigotimes_K V$ is isomorphic
to the regular representation space of $\tilde{G}$.

\bigskip

Step 2. Define

\begin{eqnarray}
&& u_0=\sum_{h\in H} x(h) \in U, \, \, v_0=\sum_{g\in G} x(g) \in V, \nonumber\\
&& \widetilde{U}=\sum_{u\in U}K \cdot u \otimes v_0 \subset U
\bigotimes_K V, \, \, \widetilde{V}=\sum_{v\in V}K \cdot u_0
\otimes v \subset U \bigotimes_K V. \nonumber
\end{eqnarray}

It is easy to see that $\widetilde{U} \oplus \widetilde{V}$ is a
faithful $\tilde{G}$-subspace of $U \bigotimes_K V$. Moreover,
when restricted to the action of $H$, the space $\widetilde{U}$ is
$H$-equivariant isomorphic to the space $U$. Similarly for
$\widetilde{V}$ and $V$ as $G$-spaces.

\medskip

Step 3. By Theorem \ref{t2.1}, $K(\widetilde{G})=K(U \bigotimes_K
V)^{\widetilde{G}}$ is rational over $K(\widetilde{U} \oplus
\widetilde{V})^{\widetilde{G}}$.

On the other hand, $K(\widetilde{U} \oplus
\widetilde{V})^{\widetilde{G}} = (K(\widetilde{U} \oplus
\widetilde{V})^H)^G$, which is $K$-isomorphic to $K(H) \cdot
K(G)$. We conclude that $K(\widetilde{G})$ is rational over $K(H)
\cdot K(G)$. (Note that the composite $K(H) \cdot K(G)$ is a free
composite, i.e. the transcendence degree of it is the sum of those
of $K(H)$ and $K(G)$.)

\medskip

Step 4. If $K(H)$ is rational (resp. stably rational) over $K$, it
is easy to see that so is $K(H) \cdot K(G)$ over $K(G)$. Thus
$K(\widetilde{G})$ is rational (resp. stably rational) over
$K(G)$.

As to the retract rationality, from the definition of retract
rationality \cite [Definition 3.2]{Sa2}, it is not difficult to
show that, (i) if $K(H)$ is retract rational over $K$, then $K(H)
\cdot K(G)$ is retract rational over $K(G)$; and (ii) if both
$K(H)$ and $K(G)$ are retract rational, then $K(H) \cdot K(G)$ is
retract rational over $K$. Hence the result.
\end{proof}

\bigskip

\bigskip

\begin{proof}[\indent Proof of Theorem \ref{t1.10}] ~

Step 1. Write $\tilde{G}=H\wr G$.

Recall the definition of the wreath product $H\wr G$.

Define $N=\bigoplus_{g\in G}H_g$ where each $H_g$ is a copy of $H$.
When we write an element $x=(\cdots, x_g, \cdots)\in N$, it is
understood that $x_g$ is the component of $x$ in $H_g$.

We will define a left action of $G$ on $N$ as follows. If $\sigma
\in G$ and $x=(\cdots, x_g, \cdots)\in N$, define ${^{\sigma}} x=
y$ where $y=(\cdots, y_g, \cdots)\in N$ with
$y_g=x_{\sigma^{-1}g}$.

The wreath product $H\wr G$ is the semi-direct product $N\rtimes
G$. More precisely, if $x,y \in N$ and $\sigma, \tau \in G$, then
$(x,\sigma) \cdot (y,\tau)=(x\cdot ({^{\sigma}} y), \sigma \tau)$.
Thus we have
\begin{eqnarray}
(\sigma x)(\tau y)=(\sigma \tau)({^{\tau^{-1}}}x \cdot y)
\end{eqnarray}
where $\sigma, \tau \in G$ and $x,y \in N$.

We will fix our notations for the group $\tilde{G}=H\wr G$, which
will be used in subsequent discussions. The groups $N$ and $G$ may
be identified (in the usual way) with subgroups of $\tilde{G}$. As
above, if $x \in N$ and $\sigma \in G$, then $(x,\sigma)$ or $x
 \sigma$ denotes an element (and the same element) in
$\tilde{G}$. For any $g \in G$, let $H_g$ be the subgroup of $N$
consisting of elements $x=(\cdots, x_{g'}, \cdots)$ satisfying the
condition that $x_{g'} =1$ for any $g'\in G \setminus \{g \}$;
define a group isomorphism $\phi_g : H \rightarrow H_g$ such that,
for any $h \in H$, if $x=\phi_g(h)$ and $x=(\cdots, x_{g'},
\cdots) \in H_g$, then $x_g=h$.

Define a subgroup $M=\sum_{g\in G \setminus \{1 \}}H_g$. Note that
the coset decomposition of $\tilde{G}$ with respect to $M$ is
given as $\tilde{G}= \cup (\sigma \cdot \phi_1(h))M$ where
$\sigma$ and $h$ run over all elements in $G$ and $H$
respectively.

\medskip

Step 2. Let $V=\bigoplus_{g\in G} K\cdot u(g)$ and
$W=\bigoplus_{x\in N} K\cdot u(x)$ be the regular representation
spaces of $G$ and $N$ respectively.

Define an action of $\tilde{G}$ on $V \bigotimes_K W$ by $(g x)
\cdot (u(g') \otimes u(y))=u(gg') \otimes u({^{g'^{-1}}}x \cdot
y)$ (following Equation (4.5)) where $g, g' \in G$ and $x,y \in
N$.

It follows that $V \bigotimes_K W$ is isomorphic to the regular
representation space of $\tilde{G}$.

\medskip

Step 3. For each $g \in G$, let $W_g=\bigoplus_{h\in H} K\cdot
u(\phi_g(h))$ be the regular representation space of $H_g$. For
any $g\in G \setminus \{1 \}$, define
\[
w_g=\sum_{h \in H}u(\phi_g(h)) \in W_g.
\]

As in Step 2 in the proof of Theorem \ref{t1.9}, we may regard
$\bigotimes_{g \in G \setminus \{1 \}} W_g$ as the regular
representation space of $M$, and regard $\bigotimes_{g\in G}W_g$
as the regular representation space of $N$, i.e. $W$. Define
\[
w'=\otimes_{g \in G \setminus \{1 \}} w_g \in \bigotimes_{g\in G
\setminus \{1 \}}W_g.
\]

Define
\[
w_0=u(1) \otimes w' \in W, \, u_0=u(1) \otimes w_0 \in V
\bigotimes_K W.
\]

Note that $x \cdot u_0=u_0$ for any $x \in M$.

\medskip

Step 4. For any $g \in G, h \in H$, define
\[
u(g;h)=(g \cdot \phi_1(h)) \cdot u_0=u(g) \otimes (u(\phi_1(h))
\otimes w') \in V\bigotimes_K W.
\]

Note that, for $g,g' \in G$ and $h,h' \in H$, we have $g \cdot
u(g';h)=u(gg';h), \phi_g(h) \cdot u(g;h')=u(g;hh')$, and
$\phi_g(h) \cdot u(g';h')=u(g';h')$ if $g \neq g'$.

For each $g \in G$, define
\[
U_g=\bigoplus_{h\in H} K\cdot u(g;h) \subset V\bigotimes_K W,
\]
and define
\[
\tilde{U}:=\bigoplus_{g\in G}U_g \subset V\bigotimes_K W.
\]

It is not difficult to show that $\tilde{U}$ is a faithful
$\tilde{G}$-subspace of $V\bigotimes_K W$. Note that $G$ permutes
the spaces $U_g$ ($g \in G$) regularly; $H_g$ acts regularly on
$U_g$, while $H_g$ acts trivially on $U_{g'}$ if $g \neq g'$.

\medskip

Step 5. Apply Theorem \ref{t2.1}. We find that $K(\tilde{G})$ is
rational over $K(\tilde{U})^{\tilde{G}}$. It remains to show that
$K(\tilde{G})$ is rational (resp. stably rational) over $K(G)$
provided that $K(H)$ is rational (resp. stably rational) over $K$.

We consider first the situation when $K(H)$ is rational over $K$.
Since $G$ permutes the spaces $U_g$ ($g \in G$) regularly, we may
choose a transcendence basis $\{v(g;i): 1 \le i \le d\}$ for
$K(U_g)^{H_g}$ (where $d$ is the order of $H$), i.e. we may write
$K(U_g)^{H_g}=K(v(g;i): 1 \le i \le d)$, such that $g \cdot
v(g';i)=v(gg';i)$ for $1 \le i \le d$.

Thus $K(\tilde{U})^{\tilde{G}}=(K(\tilde{U})^N)^G$ $=K(v(g;i):
g\in G, 1 \le i \le d)^G$. Apply Theorem \ref{t2.1}. It is easy to
see that $K(v(g;i): g\in G, 1 \le i \le d)^G$ is rational over
$K(v(g;1): g \in G)^G$, which is isomorphic to $K(G)$.

Step 6. Assume now that $K(H)$ is stably rational over $K$.
Suppose that $K(H)(w_1,$ $\cdots, w_m)$ is rational over $K$.

Define a ${\tilde{G}}$-space $\tilde{V}$ by
\[
\tilde{V}:=\bigoplus_{g\in G, 1 \le j \le m}K \cdot w(g;j).
\]
where $g \cdot w(g';j)=w(gg';j)$ and $x \cdot w(g;j)=w(g;j)$ for
any $g, g' \in G$, any $x \in N$, any $1 \le j \le m$.

Note that $\tilde{U} \bigoplus \tilde{V}$ is a faithful
${\tilde{G}}$-subspace of $(V\bigotimes_K W) \bigoplus \tilde{V}$.
By Theorem \ref{t2.1} we find that $K((V\bigotimes_K W) \bigoplus
\tilde{V})^{\tilde{G}}$ is rational over $K(V\bigotimes_K W)
^{\tilde{G}}=K({\tilde{G}})$. Again by Theorem \ref{t2.1}
$K((V\bigotimes_K W) \bigoplus \tilde{V})^{\tilde{G}}$ is rational
over $K(\tilde{U} \bigoplus \tilde{V})^{\tilde{G}}$.

Now $K(\tilde{U} \bigoplus \tilde{V})^N=\prod_{g \in G}
K(U_g)^{H_g}(w(g;j): 1 \le j \le m)$ where each $K(U_g)^{H_g}$ is
$K$-isomorphic to $K(H)$ with $g \cdot
K(U_{g'})^{H_{g'}}=K(U_{gg'})^{H_{gg'}}$ for any $g, g' \in G$. For
each $g \in G$, the field $K(U_g)^{H_g}(w(g;j): 1 \le j \le m)$ is
rational over $K$. As in Step 5, we may choose a transcendence
basis $\{v(g;i): 1 \le i \le d+m \}$ for $K(U_g)^{H_g}(w(g;j): 1
\le j \le m)$ so that $G$ acts regularly on each set $\{v(g;i): g
\in G\}$, for every $1 \le i \le d+m $. The remaining arguments are quite
similar to Step 5 and are omitted.
\end{proof}

\bigskip
\begin{prop} \label{t4.1}
Let $K$ be any field, $H\times G$ and $H\wr G$ be the direct
product and the wreath product of finite groups $H$ and $G$
respectively. If $K(H)$ is stably rational over $K$, then $K(G)$
is retract rational over $K$ if and only if so is $K(H\times G)$
{\rm(}resp. $K(H\wr G)${\rm)} over $K$.
\end{prop}

\begin{proof}
Recall a fact that, if $L_1$ and $L_2$ are stably isomorphic over
$K$, then $L_1$ is retract rational over $K$ if and only if so is
$L_2$ over $K$ \cite [Proposition 3.6]{Sa2}. Combine this fact
together with Theorem \ref{t1.9} or Theorem \ref{t1.10}.
\end{proof}

\bigskip

\begin{theorem} \label{t4.2}
Let $K$ be any field, $n$ be an odd integer, and $D_n$ be the
dihedral group of order $2n$. If $K(\bm{Z}/n\bm{Z})$ is rational
over $K$, then both $K(D_n)$ and $K(D_{2n})$ are stably rational over
$K$.
\end{theorem}

\begin{proof}
The stable rationality of $K(D_{2n})$ follows from that of
$K(D_n)$ by Theorem \ref{t1.5}.

Note that, if $n$ is an odd integer, then $(\bm{Z}/n\bm{Z})\wr
(\bm{Z}/2\bm{Z})$ is isomorphic to $(\bm{Z}/n\bm{Z})\times D_n$.
For, if $a,b \in \bm{Z}/n\bm{Z}$, $\epsilon \in
\bm{Z}/2\bm{Z}=\{\bar{0}, \bar{1} \}$ and $D_n=\langle
\sigma,\tau: \sigma^n=\tau^2=1,\
\tau\sigma\tau^{-1}=\sigma^{-1}\rangle$, the map
$\Phi:(\bm{Z}/n\bm{Z})\wr (\bm{Z}/2\bm{Z}) \rightarrow
(\bm{Z}/n\bm{Z})\times D_n$ defined by
$\Phi(2a,2b,\epsilon)=(a+b,\sigma^{a-b}\tau^{\epsilon})$ is
well-defined and is an isomorphism.

By Theorem \ref{t1.10}, the field $K((\bm{Z}/n\bm{Z})\wr
(\bm{Z}/2\bm{Z}))\simeq K((\bm{Z}/n\bm{Z})\times D_n)$ is rational
over $K$. By Theorem \ref{t1.9}, the field
$K((\bm{Z}/n\bm{Z})\times D_n)$ is rational over $K(D_n)$. Done.
\end{proof}

Remark. If $n$ is an odd integer and $K(\bm{Z}/n\bm{Z})$ is
rational over $K$, the first-named author is able to show that
$K(D_n)$ is rational over $K$ by using other methods.

\bigskip

\bigskip

\renewcommand{\refname}{\centering{References}}


\begin{thebibliography}{CHK}

\bibitem[AHK]{AHK}
H.\ Ahmad, M.\ Hajja and M.\ Kang,
\textit{Rationality of some projective linear actions},
J.\ Algebra 228 (2000), 643--658.


\bibitem[GMS]{GMS}
S.\ Garibaldi, A.\ Merkurjev and J.\ -P.\ Serre,
\textit{Cohomological invariants in Galois cohomology}, University
Lecture Series vol.\ 28, Amer.\ Math.\ Soc., Providence, 2003.

\bibitem[HK]{HK}
M.\ Hajja and M.\ Kang,
\textit{Some actions of symmetric groups},
J.\ Algebra 177 (1995), 511--535.

\bibitem[Ka]{Ka}
M.\ Kang, \textit{Actions of dihedral groups}, in ``A festschrift
in honor of Prof.\ Man-Keung Siu, 2005", Hong Kong University
Press, to appear.

\bibitem[Ku]{Ku}
H.\ Kuniyoshi, On a problem of Chevalley, Nagoya Math.\ J. 8
(1955), 65-67.


\bibitem[Pl]{Pl}
B.\ Plans, \textit{Noether's problem for $GL(2,3)$}, to appear in
``Manuscripta math.".

\bibitem[Sa1]{Sa1}
D.\ J.\ Saltman, \textit{Generic Galois extensions and problems in
field theory}, Advances in Math. 43 (1982), 250--283.

\bibitem[Sa2]{Sa2}
D.\ J.\ Saltman,
\textit{Retract rational fields and cyclic Galois extensions},
Israel J.\ Math. 47 (1984), 165--215.

\bibitem[Sw]{Sw}
R.\ G.\ Swan,
\textit{Noether's problem in Galois theory},
in ``Emmy Noether in Bryn Mawr", edited by B.\ Srinivasan and J.\ Sally,
Springer, Berlin, 1983.




\end{thebibliography}
\end{document}